\documentclass[a4paper,12pt]{article}
  \usepackage{macros_anglais}

\usepackage{amssymb}
\usepackage[english,francais]{babel}
\newtheorem{theorem}{Theorem}[section]

\newtheorem{e-proposition}[theorem]{Proposition}

\newtheorem{e-definition}[theorem]{Definition\rm}
\newtheorem{remark}{\it Remark\/}

\newenvironment{proof}{{\it Proof.} }{\hspace*{\fill}$\square$}

\newcommand{\Auttree}[1]{\Aut({\mathcal T}_#1)}

\title{On Lipschitz compactifications of trees}
\author{Beno\^{\i}t Kloeckner}

\begin{document}
\selectlanguage{english}

\maketitle

\begin{abstract}
We study the Lipschitz structures on the geodesic compactification of a 
regular tree, that are preserved by the automorphism group.
They are shown to be similar to the compactifications introduced by William Floyd,
and a complete description is given.
\end{abstract}

In \cite{Kloeckner2}, we described all possible differentiable structures
on the geodesic compactification of the hyperbolic space, for which
the action of its isometries is differentiable. 
We consider here the similar problem for regular trees and obtain
a description of ``differentiable'' compactifications, 
based on an idea of William Floyd \cite{Floyd}.
A tree has a geodesic compactification, but it is obviously not a manifold
and we shall in fact replace the differentiability condition by a Lipschitz one.

Note that we only consider regular trees so that we have a large group
of automorphisms, hence the greatest possible rigidity in our problem. A close case
is that of the universal covering of a finite graph (that is, when the automorphism
group is cocompact). Our study does not extend as it is to this case,
in particular one can
convince oneself by looking at the barycentric division of a regular tree that condition
\refeq{condition} in theorem \ref{theo} should be modified. However, 
similar results should hold, up to
considering the translates of a fundamental domain instead of the edges at some point.

This note is made of two sections. The first one recalls
some facts about regular trees and their automorphisms,
Floyd compactifications, and gives the definition of a Lipschitz
compactification. The second one contains the result and its proof.

\section{Preliminaries}

\subsection{Regular trees and their automorphisms}

We denote by $\mathcal{T}_n$ the regular tree of valency $n\geq 3$ and by
$T_n$ is topological realization, obtained by replacing each abstract edge by a segment.
All considered metrics on $T_n$ shall be length metrics, since
general metrics could have no relation at all with the combinatorial structure of
$\mathcal{T}_n$. Up to isometry, two length metrics on $T_n$ that are compatible
with the topology differ only by the length of the edges. We shall therefore
identify $T_n$ equipped with such a metric and ${\mathcal T}_n$ equipped with
a labelling of the edges by positive real numbers (the label corresponding to the
length of the edge). When all edges are chosen
of length $1$, we call the resulting metric space the ``standard metric realization''
of ${\mathcal T}_n$, denoted by $T_n(1)$. Its metric shall be denoted by $d$;
it coincides on vertices with the usual combinatorial distance.

There is a natural one-to-one correspondence between automorphism of
${\mathcal T}_n$ and isometries of $T_n(1)$. We denote both groups
by $\Auttree{n}$ and endow them with the compact-open topology, so that a basis 
of neighborhoods of identity is given by the sets
$B_K(\id)=\ensemble{\phi\in\Auttree{n}}{\phi(x)=x\quad\forall x\in K}$
where $K$ runs over all finite sets of vertices.

Given an automorphism $\phi$,
one defines the \emph{translation length} of $\phi$ as the integer
$T(\phi)=\min_x\{d(x,\phi(x))\}$
where the minimum is taken over all points (not only vertices) of $T_n(1)$.
The following alternative is classical:
\begin{enumerate}
\item if $T(\phi)>0$ then there is a unique invariant bi-infinite path
      $(x_i)_{i\in\mZ}$ and $\phi(x_i)=x_{i+T(\phi)}$ for all $i$,
\item if $T(\phi)=0$ then either $\phi$ fixes some vertex, or
      $\phi$ has a unique fixed point in $T_n(1)$, which  is the midpoint
      of an edge.
\end{enumerate}
In the first case, $\phi$ is said to be a \emph{translation} (a 
\emph{unitary} translation if $T(\phi)=1$).  Any translation
is a power of a unitary translation.

\subsection{Compactification of trees}

The standard metric tree $T_n(1)$ is a CAT($0$) complete length space, thus
is a Hadamard space (see for example \cite{Burago}). Therefore, it has a geodesic 
compactification we now briefly describe.

A boundary point $p$ is a class of asymptotic
geodesic rays, where two geodesic rays $\gamma_1=x_0,x_1,\dots,x_i,\dots$ and 
$\gamma_2=y_0,y_1,\dots,y_j,\dots$ are said to be \emph{asymptotic}
if they are eventually identical: there are indices
$i_0$ and $j_0$ so that for all $k\in\mN$, on has
$x_{i_0+k}=y_{j_0+k}$. The point $p$ is said to be the endpoint
of any geodesic ray of the given asymptoty class.

The union $\overline{T}_n=T_n\cup\partial T_n$ is given the following topology:
for a point that is not on the boundary, a basis of neighborhoods is given by its neighborhoods in 
$T_n$; for a boundary point $p$, a basis of neighborhoods is given by the 
connected components
of $T_n\setminus\{x\}$ containing a geodesic ray asymptotic to $p$, 
where $x$ runs over the vertices 

It is a general property of Hadamard spaces that $\Auttree{n}$ acts on
$\overline{T}_n$ by homeomorphisms for this topology.
Our goal will be to see which additional structure can be added to this topology,
that is preserved by $\Auttree{n}$.

We have no differentiable structure on $\overline{T}_n$, but due to
the Rademacher theorem it is natural to look at Lipschitz structures instead.

\begin{e-definition}
Let $X$ be a metrizable topological space. A Lipschitz structure $[\delta]$ on $X$
is the data of a metric $\delta$ that is compatible with the topology of $X$,
up to local Lipschitz equivalence (two metrics $\delta_1$, $\delta_2$ are said to be locally
Lipschitz equivalent if the identity map $(X,\delta_1)\to (X,\delta_2)$ is locally bilipschitz).
\end{e-definition}

The natural isomorphisms of a space $X$ endowed with a Lipschitz structure are
the locally bilipschitz maps.
Usually, for an action of a Lie group on a manifold to be differentiable, one asks
the map $G\times M\to M$ to be differentiable. Similarly, we say that an action
of a topological group $\Gamma$ on a metrizable topological space $X$ is Lipschitz
if it is a continuous action by locally bilipschitz maps, and if moreover
the Lipschitz factor is locally uniform.

We can now define our main object of study.

\begin{e-definition}
A Lipschitz compactification of ${\mathcal T}_n$ is a Lipschitz structure
$[\delta]$ on $\overline{T}_n$, where $\delta$ is a length metric,
and such that the action of $\Auttree{n}$ on $\overline{T}_n$ is Lipschitz.
\end{e-definition}


In \cite{Floyd}, Floyd introduced a method for compactifying a graph. We give
definitions that are adapted to the simpler case of trees.

\begin{e-definition}
By a Floyd function we mean a function $h:\mN\to ]0,+\infty[$ such that 
$\sum_r h(r) <+\infty$.
Two Floyd functions $h_1$, $h_2$ are said to be \emph{comparable} if there is a $C>1$ such that
for all $r\in\mN$ one has $C^{-1} h_2(r)\le h_1(r) \le C h_2(r)$.
\end{e-definition}

\begin{e-definition}
A Floyd metric on $\overline{T}_n$ is the length metric obtained from a vertex $x_0$ and a Floyd function $h$
by assigning to each edge $e$ the length $h(d)$, where $d\in\mN$ is the combinatorial distance between
$e$ and $x_0$.

By a Floyd compactification of $\mathcal{T}_n$ we mean the topological space $\overline{T}_n$
endowed with the Lipschitz structure corresponding to a Floyd metric.
\end{e-definition}

The condition that $\sum h(r)$ converges ensures that we do get a distance on $\overline{T}_n$. For example,
the distance
between two boundary points $p$ and $p'$ is $2\sum_{r\ge R}h(r)$ where $R$ is the combinatorial
distance between $x_0$ and the only geodesic joining $p$ and $p'$.

Two Floyd metrics obtained from the same point $x_0$ and Floyd functions
$h_1$, $h_2$ are easily seen to define the same Lipschitz structures
if and only if $h_1$ and $h_2$ are comparable.

\section{Description of all Lipschitz compactifications of regular trees}

\begin{theorem}\label{theo}
Any Lipschitz compactification of $\mathcal{T}_n$ is a Floyd compactification.

The Floyd compactification of $\mathcal{T}_n$ obtained from a Floyd function
$h$ and a base point $x_0$ is a Lipschitz compactification if and only if
there is a constant $0<\eta<1$ so that for all $r\in\mN$
\begin{equation}
h(r+1)\ge\eta\, h(r).
\label{condition}
\end{equation}
\end{theorem}

\begin{remark}
Condition \refeq{condition} implies that $h$ decreases at most exponentially fast. It is interesting
to compare this with the usual conformal compactification of the hyperbolic space, obtained
by multiplying the metric by a factor that is exponential in the distance to a fixed point.
\end{remark}
\begin{remark}
Condition \refeq{condition} implies that the considered Lipschitz structure depends only upon $h$, 
not $x_0$. We can therefore denote this compactification by $\overline{T}_n(h)$.
\end{remark}

\begin{proof}
We first prove that any Lipschitz compactification of $\mathcal{T}_n$ is a Floyd compactification.

Let $\delta'$ be any length metric in the given Lipschitz class, and fix any vertex $x_0$ of $\mathcal{T}_n$.
We define $h$ by $h(r)=\min\delta'(x,y)$ where the minimum is taken over all edges $xy$ that are at
combinatorial distance $r$ from $x_0$. Then $h$ is a Floyd function because $x_0$ is at finite $\delta'$
distance from the boundary. Denote by $\delta$ the Floyd metric obtained
from $x_0$ and $h$, and let us prove that $[\delta]=[\delta']$. It is sufficient to
prove that there is a constant $C$ so that for all $r$, two edges that are at 
combinatorial distance $r$ from $x_0$ have their $\delta'$ lengths that differ by a factor at most $C$.

For any $R\in\mN$, let $B(R)$ be the closed ball of radius $R$ and center $x_0$ in $T_n(1)$. It contains a 
finite number of edges, so that there is a constant $C_R$ that satisfies the above property
 for all $r\le R$.

Since the compactification is assumed to be Lipschitz, for all $p\in\partial T_n$ there are
a neighborhood $V$ of $p$, a neighborhood $U$ of the identity and a constant $k$ so that
any $\phi\in U$ is $k$-Lipschitz on $V$. Since $\partial T_n$ is compact, we can find a finite
number of such quadruples $(p_i,V_i,U_i,k_i)$ so that the $V_i$ cover $\partial T_n$. Moreover
we can assume that the $V_i$ are the connected components of $T_n\setminus B(R)$ for some
radius $R$, and that $U=\cap U_i=B_{B(R)}(\id)$. Since for all $i$ and $r>R$, $U$ acts transitively
on the set of edges of $V_i$ that are at combinatorial distance $r$ from $x_0$, those edges
have their $\delta'$-length that differ by a factor at most $C'=\sup k_i$. Moreover,
there is an automorphism $\phi_0$ that fixes $x_0$ and permutes cyclically the $V_i$.
Since $\phi_0$ is locally Lipschitz, there is a $R'$ and a $C''$ so that for all $r\ge R'$
and all couple $(i_1,i_2)$, there are edges of $V_{i_1}$ and $V_{i_2}$ that are
at combinatorial distance $r$ from $x_0$ and whose $\delta'$ lengths differ by a
factor at most $C''$. The supremum $C$ of $C_{R'}$ and $C''C'^2$ is the needed constant.

Consider now the Floyd compactification obtained from $x_0$ and $h$ and denote by $\delta$ 
the associated Floyd metric.
By construction, any automorphism $\phi$ of $\mathcal{T}_n$ that fixes $x_0$
is an isometry for $\delta$, thus
is locally bilipschitz for the corresponding Lipschitz structure.

Two translations are close to one another when they differ by an element
close to identity. An element close enough to identity must fix
$x_0$, thus is an isometry. Therefore, we only need to prove that a given translation
is Lipschitz to get that all automorphisms in a neighborhood are equilipschitz.
Checking unitary translations is sufficient since any translation is an iterate
of one of those.

Let $\phi$ be a unitary translation, and $\gamma=\dots,y_{-1},y_0,y_1,\dots$ 
be its translated geodesic, where we assume that $y_0$ realizes the minimal
combinatorial distance $d_0$ between vertices of $\gamma$ and $x_0$. 
By local finiteness, $\phi$ is locally bilipschitz around
any point of $T_n$ and we need only check the boundary.

Let us start with the attractive endpoint $p$ of $\gamma$. Assume that
our Floyd compactification is Lipschitz. It implies that $\phi$ is locally
bilipschitz around $p$, in particular there is a $r_0>0$ and a $k>1$ such that
for any $r\ge r_0$,
\begin{eqnarray*}
k\,\delta(y_{r+1},y_{r+2}) &\ge& \delta(y_r,y_{r+1}) \\
h(r+d_0+1) &\ge& k^{-1}h(r+d_0)
\end{eqnarray*}
which gives condition \refeq{condition}.

Conversely, assume that \refeq{condition} holds.

For any vertex $x$ we have
\[\left|d(\phi(x),x_0)-d(x,x_0)\right|\le 1+2 d_0\]
since the worst case is when $x=x_0$ or $x$ is in a connected component of
$\mathcal{T}_n\setminus\{x_0\}$ other than that of $\gamma$.
Therefore, the length of an edge and of its image by $\phi$ differ by a factor
bounded by $\eta^{-(1+2d_0)}$. Therefore, $\phi$ is Lipschitz. Since $\phi^{-1}$ is also a unitary translation,
$\phi$ is bilipschitz.
\end{proof}

It would be interesting to consider more general spaces, for example euclidean buildings
or CAT(-1) buildings like the $I_{pq}$ described
by Bourdon in \cite{Bourdon}. It is not obvious how to define the Floyd compactification:
for example, a mere scaling of the distance in each cell by a factor depending on the combinatorial
distance to a fixed cell would create gluing problems (an edge shared by two faces having
two different length). This spaces could therefore be less flexible than trees.

\bibliographystyle{plain}
\bibliography{biblio.bib}

\end{document}